\documentclass{article}

\usepackage{amsmath, amsfonts, amssymb, amsthm}
\usepackage[small]{titlesec}

\newcommand{\nablabar}{\overline{\nabla}}

\newcommand{\qbar}{\overline{q}}
\newcommand{\jbar}{\overline{j}}

\newcommand{\lbar}{\overline{l}}

\newcommand{\ov}[1]{\overline{#1}}
\newcommand{\ghat}{\hat{g}}
\newcommand{\Ric}{\mathrm{Ric}}
\newcommand{\Rm}{\mathrm{Rm}}
\newcommand{\Stilde}{\tilde{S}}
\newcommand{\Btilde}{\tilde{B}}
\newcommand{\ddt}{\frac{\partial}{\partial t}}

\newcommand{\myle}{\phantom{l}\le\phantom{l}}

\newtheorem{theorem}{Theorem}[section]

\newtheorem{corollary}[theorem]{Corollary}

\numberwithin{equation}{section}

\begin{document}

\centerline{\bf  \large Interior derivative estimates for the K\"ahler-Ricci flow\footnote{Research supported in part by NSF grants DMS-0848193 and  DMS-1105373.   In addition, the first-named author was supported in part by a California State Faculty Support Grant (SFSG) and the second-named author  by a Sloan Fellowship.  }}

\bigskip

\centerline{Morgan Sherman$^*$ and Ben Weinkove$^\dagger$}

\begin{abstract}
We give a maximum principle proof of  interior derivative estimates for  the K\"ahler-Ricci flow, assuming local uniform bounds on the metric.
\end{abstract}

\section{Introduction}

Let $(M, \hat{\omega})$ be a K\"ahler manifold of complex dimension $n$.  Let $\omega=\omega(t)$ be a solution of the K\"ahler-Ricci flow on $M \times [0,T]$, for some $T>0$:
\begin{equation} \label{krf}
\ddt \omega = - \Ric(\omega), \qquad \omega|_{t=0} = \omega_0,
\end{equation}
with $\omega_0$ a smooth initial K\"ahler metric.

Fix a point $p \in M$ and denote by  $B_r \subset M$ the open ball centered at $p$ of radius $r$ for $0<r<1$ with respect to $\hat{\omega}$.  We assume that $r$ is sufficiently small so that $\overline{B_r}$ is contained in a single holomorphic coordinate chart.  Our main result is as follows:

\begin{theorem} \label{maintheorem} Let $N>1$ satisfy
\begin{equation} \label{assume}
\frac{1}{N} \hat{\omega} \le \omega \le N \hat{\omega}, \quad \textrm{on } \ov{B_r} \times [0,T].
\end{equation}
Then for each $m=0, 1, 2, \ldots$ there exist constants $C$ and $C_m$ depending only on $\hat{\omega}$ and $T$
such that on $B_{r/2} \times (0,T]$,
\begin{enumerate}
\item[(i)] $\displaystyle{| \hat{\nabla} \omega|_{\omega}^2 \le C \frac{N^3}{r^2 t}}$, for $\hat{\nabla}$ the covariant derivative with respect to $\hat{\omega}$.
\item[(ii)] $\displaystyle{|\Rm|_{\omega}^2 \le C_0  \frac{N^8}{r^4 t^2} }$.
\item[(iii)] $\displaystyle{|\nabla_{\mathbb{R}}^m \Rm|_{\omega}^2 \le C_m \left( \frac{N^4}{r^2 t} \right)^{m+2}}$ for $m=1, 2, \ldots$, for $\nabla_{\mathbb{R}}$ the real covariant derivative with respect to the metric $\omega$.
\end{enumerate}

Moreover, if we allow the constants $C$ and $C_m$ to depend also on $\omega_0$ then the estimates (i), (ii) and (iii) hold with each factor of $t$ on the right hand side  replaced by $1$.
\end{theorem}

We prove this result using the maximum principle.  Note that by work of Shi \cite{Shi, Shi2} it was already known that a bound on curvature as in (ii)  implies (iii) (nevertheless, we include a proof here, for the sake of completeness).  Theorem \ref{maintheorem} implies the following:

\begin{corollary} \label{corollary}
Let $N>1$ satisfy
\begin{equation} \label{assume2}
\frac{1}{N} \hat{\omega} \le \omega \le N \hat{\omega}, \quad \textrm{on } \ov{B_r} \times [0,T].
\end{equation}
Then for each $m=0, 1, 2, \ldots$ there exist constants $C_m$, $\alpha_m$, $\beta_m$ and $\gamma_m$  depending only on $m$, $\hat{\omega}$ and $T$ such that 
\begin{equation} \label{mainestimate}
| \hat{\nabla}_{\mathbb{R}}^m \omega |_{\hat{\omega}} \le C_m \frac{N^{\alpha_m}}{r^{\beta_m} t^{\gamma_m}}, \quad \textrm{on } B_{r/2} \times (0,T],
\end{equation}
Moreover, if we allow the constants $C_m$, $\alpha_m$ and $\beta_m$ to depend also on $\omega_0$ then (\ref{mainestimate}) holds with $\gamma_m=0$.
\end{corollary}

Namely, a local uniform estimate for the metric along the K\"ahler-Ricci flow implies local derivative estimates to all orders.  This fact in itself is not new.  Indeed the local PDE theory of Evans-Krylov  \cite{E, K} can be applied to the K\"ahler-Ricci flow equation (see for example \cite{Ch} or the generalization in  \cite{G}).   The key point here is to establish this via Theorem \ref{maintheorem} whose proof uses only elementary maximum principle arguments.

The form of the estimate (\ref{mainestimate}) may be useful for applications and does not seem to be written down explicitly elsewhere in the literature.  
 When considering the K\"ahler-Ricci flow on projective varieties, it is often the case that one obtains a uniform estimate for the metric $\omega$ away from a subvariety (see for example \cite{ST, SW0, SW, SW2, TZ, Ts, Z}).   Theorem \ref{maintheorem} can be used to replace global arguments.    To illustrate, suppose that $\omega=\omega(t)$ solves the K\"ahler-Ricci flow on a compact K\"ahler manifold $M$ and there exists an analytic hypersurface $D \subset M$ whose associated line bundle $[D]$ admits a holomorphic section $s$ vanishing to order 1 along $D$. Assume that 
 \begin{equation}
\frac{1}{C} |s|_H^{\alpha} \hat{\omega} \le \omega \le \frac{C}{ |s|^{\alpha}_H} \hat{\omega}, \quad \textrm{on } (M \setminus D) \times [0,T]
\end{equation}
for some positive constants $C$ and $\alpha$, where $H$ is a Hermitian metric on $[D]$.  An elementary argument shows that Theorem \ref{maintheorem} implies the existence of $C_m$, $\alpha_m$ and $\gamma_m$ such that
 \begin{equation} \label{De}
| \hat{\nabla}_{\mathbb{R}}^m \omega |_{\hat{\omega}}\le \frac{C_m}{ t^{\gamma_m} |s|^{\alpha_m}_H}, \quad \textrm{on } (M \setminus D) \times (0,T]
\end{equation}
for each $m=1, 2, \ldots$.  Moreover we can take $\gamma_m=0$ if we allow $C_m$ and $\alpha_m$ to depend on the initial metric $\omega_0$.  Estimates of the form of (\ref{De}) are used for example in \cite{SW, SW2}.  In particular, Corollary \ref{corollary} gives an alternative proof of the results in Section 4 of \cite{SW}.

Finally we remark that since our result is completely local, we may and do assume that $M =\mathbb{C}^n$, $p=0$ and $\hat{\omega}$ is the Euclidean metric.
 We will write $g$ and $\hat{g}$ for the K\"ahler metrics associated to $\omega$ and $\hat{\omega}$.
All magnitudes $|\cdot|$ are taken with respect to the metric $g$.  We shall use the letter $C$ (as well as $C', C'',$ etc.) for a uniform constant (depending only on $m$, $\hat{\omega}$, and $T$) which may differ from line to line.

In Sections \ref{secfirst}, \ref{seccurv} and \ref{seccurv2} we prove parts (i), (ii) and (iii) of Theorem \ref{maintheorem} respectively.   In Section \ref{secproof} we give a proof of Corollary \ref{corollary}.

\section{Bound on the first derivative of the metric} \label{secfirst}

In this section we prove the estimate on the first derivative of the metric $g$, establishing part (i) of Theorem \ref{maintheorem}.  This gives a local parabolic version of the well-known Calabi `3rd order' estimate \cite{Ca} for the complex Monge-Amp\`ere equation  (used by Yau \cite{Y} in his solution of the Calabi conjecture).  There exist now many generalizations of Calabi's estimate  (see for example \cite{Che, To, TWY,   ZZ}).  A global parabolic Calabi estimate was applied to the case of the K\"ahler-Ricci flow in \cite{Cao}.  Phong-Sesum-Sturm \cite{PSS} later gave a neat and explicit computation in this  which we will make use of here for our local estimate.

We wish to bound the quantity 
\begin{equation}
S = | \hat{\nabla} g|^2 = g^{i \ov{j}} g^{k \ov{l}} g^{p \ov{q}} \hat{\nabla}_i g_{k \ov{q}} \ov{ \hat{\nabla}_j g_{l \ov{p}}}
\end{equation}
 where we write $\hat{\nabla}$ for the covariant derivative with respect to $\hat{g}$.
Write $r_0=r$ and let $\psi$ be a nonnegative $C^{\infty}$ cut-off function that is identically equal to 1 on $\ov{B_{r_1}}$ and vanishes outside $B_r$, where $r_0 >r_1>r/2$.  We may assume that 
\begin{equation}
	|\nabla\psi|^2, \, |\Delta\psi| \le C \frac{N}{r^2},
\label{bound for grad phi}
\end{equation}
where $\Delta = \nabla^{\ov{j}} \nabla_{\ov{j}}  = g^{p \ov{q}} \nabla_p \nabla_{\ov{q}}$.
Thus
\begin{equation} 
	(\partial_t-\Delta)(\psi^2 S) \le 
		\psi^2 (\partial_t-\Delta) S + C \frac{N}{r^2} S 
		+ 2 \left| \langle \nabla \psi^2,\, \nabla S \rangle \right|,
\label{preliminary inequality on evol of psi squared S}
\end{equation}
where we are writing $\langle \nabla F, \nabla G \rangle = g^{i \ov{j}} \partial_i F \partial_{\ov{j}} G$ for functions $F, G$.
Following the notation in \cite{PSS}, we introduce the endomorphism 
$h^{i}{}_{k} = \hat{g}^{i \ov{j}} g_{\ov{j} k}$ and let $X$ be the tensor with components $X_{il}^{k} = (\nabla_{i} h \cdot h^{-1})^{k}{}_{l}$, so that $S = |X|^2$.   Note that $X$ is the difference of the Christoffel symbols of $g$ and $\hat{g}$.

  An application of Young's inequality gives 
\begin{equation}
	2 \left| \langle \nabla \psi^2,\, \nabla S \rangle \right|
		\le \psi^2 ( |\nabla X|^2+|\nablabar X|^2 ) 
		+ C \frac{N}{r^2} S .
\label{cross term estimate in evolution of psi squared S}
\end{equation}
We now use the evolution equation for $S$ derived by Phong-Sesum-Sturm (see equation (2.51) of \cite{PSS}) which, in the case where $\hat{\omega}$ is Euclidean, has the simple form:
\begin{equation} 
	(\partial_t-\Delta)S \ =\ -\left( |\nabla X|^2 + |\nablabar X|^2 \right).
\label{evolution of S}
\end{equation}
Combining 
(\ref{preliminary inequality on evol of psi squared S}, 
\ref{cross term estimate in evolution of psi squared S}, 
\ref{evolution of S}) 
we find 
\begin{equation} 
	(\partial_t-\Delta)(\psi^2 S) 
	\le C \frac{N}{r^2}S .
\label{combined est}
\end{equation}

We now need to use the evolution equation for $\mathrm{tr}\, h$ from  \cite{Cao}, which is a parabolic version of an estimate from \cite{A, Y}.   More precisely, we can apply equations (2.28) and (2.31) of \cite{PSS} and use the fact that the fixed metric is Euclidean to obtain
\begin{equation} \label{Delta-dt Delta hat phi}
	(\partial_t - \Delta)(\mathrm{tr}\,h) 
	= -\ghat^{i \ov{j}} g^{k  \ov{l}} g^{p \ov{q}} 
		\hat\nabla_i g_{\ov{l} p} \overline{ \hat\nabla_j g_{\ov{k} q} }.
\end{equation}
Hence
\begin{equation} \label{dt-Delta Delta hat phi est}
	(\partial_t-\Delta)(\mathrm{tr}\,h) \ \le\ - \frac{S}{N} .
\end{equation}
Let $f(t)$ denote either the function $t$ or the constant $1$.  Then $0 \le f(t) \le \max(T,1)$ and $f'(t) = 1$ or $0$ so that we get, for any positive constant $B$, 
\begin{equation*} 
	(\partial_t-\Delta)(f(t) \psi^2 S + B \, \mathrm{tr}\,h) 
	\ \le\ C \frac{N}{r^2}S - \frac{B}{N}S .
\end{equation*}
Let $B = \frac{N^2}{r^2}(C+1)$.  Then by the maximum principle,   the maximum of $f(t) \psi^2 S + B \, \mathrm{tr}\,h$ on $\ov{B_r} \times [0,T]$ can only occur at $t=0$ or on the boundary of $\ov{B_r}$, where $\psi=0$.  Since $\textrm{tr} \, {h} \le nN$, we have
\begin{equation} 
	S \le C \frac{N^3}{f(t) r^2} \ \ \text{on}\ \ov{B_{r_1}} \times (0,T].
\label{bound for S on B_r1}
\end{equation}
giving part (i) of Theorem \ref{maintheorem}.

\section{Bound on curvature} \label{seccurv}

We now prove part (ii) of Theorem \ref{maintheorem}.  For global estimates of this type, see for example \cite{Chau, PSSW}.  We fix a smaller radius $r_2$ satisfying $r_1 > r_2 > r/2$.  In this section we let $\psi$ be a cut-off function, identically 1 on $\ov{B_{r_2}}$ and identically 0 outside $B_{r_1}$.  As before we may assume $|\Delta \psi|, |\nabla\psi|^2 \le C N /r^2$ for some uniform constant $C$.  Calculate
\begin{align}
	(\partial_t-\Delta) R_{\jbar i \lbar k} = &
		- R_{\jbar i}{}^{p \qbar} R_{\lbar k \qbar p} 
		+ R_{\lbar i}{}^{p \qbar} R_{\jbar k \qbar p} 
		- R_{\jbar p \lbar}{}^{\qbar} R^{p}{}_{i \qbar k} \nonumber\\
		& - R_{\jbar p} R^p{}_{i \lbar k} - R_{\lbar p} R_{\jbar i}{}^p{}_k, 
\label{laplacian of curvature}
\end{align}
and therefore (cf. \cite{H})
\begin{align}
	(\partial_t-\Delta) |\Rm|^2 \le  
 - | \nabla \Rm |^2 - | \nablabar \Rm |^2 + C | \Rm|^3,
\end{align}
where we are writing $| \Rm|^2 = R_{\ov{j} i \ov{l} k} R^{i \ov{j} k \ov{l}}$ etc.

As before we set $f(t) = t,1$.  We introduce the function 
\begin{equation}
	\Stilde = f S + C_1 N \,\mathrm{tr}\,h 
\label{Stilde definiton}
\end{equation}
where $C_1$ is a large uniform constant.  Note that by (\ref{bound for S on B_r1}) we have $\Stilde \le C \frac{N^3}{r^2}$ at every $(x,t) \in \ov{B_{r_1}} \times [0,T]$.  Furthermore $\Stilde$ satisfies 
\begin{equation}
(\partial_t-\Delta)\Stilde \le -f (|\nabla X|^2+|\nablabar X|^2) - C_2 S
\end{equation}
where $C_2 = C_1-f' \gg 1$ is uniform.  Let $K = C_3 N^4/r^2$ where $C_3 \gg 1$ is a uniform constant.  Note that we may assume $K/2 \le K - \Stilde \le K$.  We will establish our bound for $|\Rm|$ by using a maximum principle argument for the function $F = f^2 \frac{\psi^2|\Rm|^2}{K-\Stilde} + \Btilde\Stilde$ where $\Btilde = C_4/N^3$ with $C_4 \gg 1$ uniform.  We begin by computing 
\begin{align}
	(\partial_t-\Delta) \left( \psi^2 \frac{|\Rm|^2}{K-\Stilde} \right)
		=& - \Delta \psi^2 \frac{|\Rm|^2}{K-\Stilde} 
			+ \psi^2 \frac{(\partial_t-\Delta)|\Rm|^2}{K-\Stilde} 
			+ \psi^2 \frac{(\partial_t-\Delta)\Stilde}{(K-\Stilde)^2} |\Rm|^2 
			\nonumber\\
		& - 2\psi^2 \frac{|\nabla\Stilde|^2|\Rm|^2}{(K-\Stilde)^3}
			- 4\mathrm{Re} \frac{\psi \langle \nabla\psi , \nabla |\Rm|^2\rangle}{K-\Stilde}
			\nonumber\\
		& - 4\mathrm{Re} \frac{\psi \langle \nabla\psi,  \nabla \Stilde \rangle |\Rm|^2}{(K-\Stilde)^2}
			-2 \mathrm{Re} \frac{\psi^2\langle \nabla |\Rm|^2,  \nabla \Stilde \rangle}{(K-\Stilde)^2}
\label{evolution of complicated fn}
\end{align}
and thus
\begin{align} \nonumber
\lefteqn{	(\partial_t-\Delta)  \left( \psi^2 \frac{|\Rm|^2}{K-\Stilde} \right)} \\
		& \myle \frac1{(K-\Stilde)^2} \Biggl[ |\Delta\psi^2| (K-\Stilde) |\Rm|^2
		 + \psi^2(K-\Stilde) \left( C|\Rm|^3 - |\nabla\Rm|^2 - |\nablabar\Rm|^2 \right) 
			\nonumber\\
		& + \psi^2 \left( -f|\nabla X|^2 -f|\nablabar X|^2 - C_2 S \right) |\Rm|^2 
			- 2 \psi^2 \frac{|\nabla\Stilde|^2|\Rm|^2}{K-\Stilde} 
			\nonumber\\
		& + 16 |\nabla\psi|^2(K-\Stilde)|\Rm|^2 
			+ \frac12\psi^2(K-\Stilde)\left|\nabla \Rm \right|^2 + \frac12\psi^2(K-\Stilde)\left| \ov{\nabla} \Rm \right|^2
			\nonumber\\
		& + \frac1{K-\Stilde}\psi^2|\nabla\Stilde|^2|\Rm|^2 
			+ 4 (K-\Stilde) |\nabla\psi|^2|\Rm|^2
			\nonumber\\
		& + \frac{4}{K-\Stilde} \psi^2 |\nabla\Stilde|^2|\Rm|^2
			+ \frac{1}{2} \psi^2 (K-\Stilde) \left|\nabla \Rm\right|^2 +\frac{1}{2} \psi^2 (K-\Stilde) \left|\nablabar \Rm\right|^2
			\Biggr] .
\label{evolution of complicated fn, bound 1}
\end{align}
We wish to bound (\ref{evolution of complicated fn, bound 1}) in terms of $|\Rm|^2$.  Label the terms $(1), (2), \ldots, (16)$.  Then the bad terms are $(1), (2)$ and (9) through (16) while the remaining terms are all good.  One sees that $(1) + (9) + (13) \le C \frac{N}{Kr^2} |\Rm|^2$ while 
$[(10)+(11) + (15)+(16)] + [(3)+(4)] \le 0$ and $(12) + \frac12 (8)\le 0$.  It remains only to bound the terms $(2)$ and $(14)$.  For $(2)$ we argue as follows:  we may assume that at a maximum for the function $F$ we have a lower bound of the form 
\begin{equation}
	f |\Rm| \ge C K, \ C \gg 1 
\label{assumed minimum for f Rm}
\end{equation}
for if not we can apply a maximum principle argument immediately:  At any $(x,t) \in \ov{B_{r_1}} \times (0,T]$ we would have $F$ is at most $C K + C/r^2$ which implies that 
$$ f^2 |\Rm|^2 \le C \frac{N^8}{r^4} \textrm{ on } \ov{B_{r_2}} \times (0,T] . $$
Now since $\hat{\omega}$ is Euclidean we have
\begin{equation}
	\left| \nablabar X \right|^2= |\Rm - \widehat{\Rm}|^2  = | \Rm|^2.
\label{bound for nabla bar X in terms of Rm}
\end{equation}
Hence, using (\ref{assumed minimum for f Rm}), we have $(2) + \frac12 (6) \le 0$.  Finally, to control $(14)$ we use 
\begin{equation}
| \nabla\Stilde |^2 \le 4 f^2 S ( |\nabla X|^2 + |\nablabar X|^2 ) + 2n C_1^2 N^4 S .
\end{equation}
Here we have made use of a well-known estimate (computed in \cite{Y}) which implies that 
$|\nabla\,\mathrm{tr}\,h|^2 \le n N^2 S$.  Now we find $(14) + \frac12 [(5) + (6) + (7)] \le 0$ if in $K = C_3 N^4/r^2$ we choose $C_3 \gg C_1$.  In total then we have
\begin{equation}
	(\partial_t-\Delta) \left( \frac{\psi^2 |\Rm|^2}{K-\Stilde} \right) 
	\le \frac{C}{N^3} |\Rm|^2 .
\end{equation}
Therefore
\begin{equation}
	(\partial_t-\Delta) \left(  \frac{\psi^2 f^2 |\Rm|^2}{K-\Stilde} 
		+ \tilde{B} \Stilde \right) 
		\ \le\ - \frac{f}{N^3} |\Rm|^2,
\label{bound on final complicated fn}
\end{equation}
if in $\tilde{B} = C_4/N^3$ we pick $C_4$ large enough.  This implies that the maximum of $F$ on $\ov{B_{r_1}} \times [0,T]$ can only occur at $t=0$ or on the boundary of $\ov{B_r}$, where $\psi=0$.  Hence $F$ is bounded above by $C/r^2$.
  Therefore at any $(x,t)$ in $\ov{B_{r_2}}\times [0,T]$ we have 
$f^2 |\Rm|^2 \le   C' N^4 / r^4$.  Comparing with our comments following (\ref{assumed minimum for f Rm}) we arrive at the following estimate:
\begin{equation}
	|\Rm|^2 \le C \frac{N^8}{f(t)^2 r^4}\ \ \text{on}\  \ov{B_{r_2}} \times (0,T].
\label{bound on curvature}
\end{equation}

\section{Higher order estimates} \label{seccurv2}

We finish the proof of Theorem \ref{maintheorem} by establishing  bounds on the derivatives of curvature,  following the basic idea of Shi  \cite{Shi, Shi2} (cf.  \cite{B, CK, CLN}).  Our setting here is slightly different from that of Shi, where it is assumed that curvature is uniformly bounded (independent of $t$) but that (\ref{assume}) does not necessarily hold.  Although the result we need can be recovered from what is known in the literature, we include the short proof for the sake of completeness.  
    Fix a sequence of radii $r = r_0 > r_1 > r_2 > \ldots > r/2$.  For a fixed $m$ we will denote by $\psi$ a cutoff function which is zero outside $B_{r_{m+1}}$ and identically 1 on $\ov{B_{r_{m+2}}}$.  

We now work in real coordinates, writing, in this section, $\nabla$ for the real covariant derivative $\nabla_{\mathbb{R}}$.  Write $\nabla^m$ for $\nabla \nabla \cdots \nabla$ ($m$ times). The key evolution equation we need is due to Hamilton \cite{H}:
\begin{equation} 
	(\partial_t-\Delta)| \nabla^m\Rm |^2 
	= - |\nabla^{m+1}\Rm|^2 + \sum_{i+j=m} \nabla^i \Rm * \nabla^j \Rm * \nabla^m \Rm,
\label{Higher order evolution general formula}
\end{equation}
where we are writing $S*T$ to denote a linear combination of the tensors $S$ and $T$ contracted with respect to the metric $g$.  To clarify (\ref{Higher order evolution general formula}), we take $\Delta$ here to be the complex Laplacian, which, acting on functions, is half the usual Riemannian Laplace operator.  When comparing to the formula in \cite{H} note that Hamilton's Ricci flow equation includes a factor of 2 which is not present in our equation (\ref{krf}).

We will show inductively that 
\begin{equation}
	\left| \nabla^m \Rm \right|^2 \le C \left( \frac{N^4}{f(t) r^2} \right)^{m+2}
	\quad \textrm{on } \ov{B_{r_{m+2}}} \times (0,T]
\label{Bound on grad^m Rm}
\end{equation}
for every $m \ge 0$, the base case $m=0$ having already been established in Section 
\ref{seccurv}.  Assume (\ref{Bound on grad^m Rm}) holds for every value $< m$.  Let 
$A = N^4/r^2$.
We will apply the maximum principle argument to the function 
\begin{equation}
	F = \psi^2 f^{m+2} | \nabla^m\Rm |^2 + B f^{m+1} | \nabla^{m-1}\Rm |^2 
\end{equation}
where $B = C_1 A$ with $C_1 \gg 1$ a large uniform constant.  Let $(x_0,t_0) \in \ov{B_{r_{m+1}}} \times [0,T]$ be the point at which $F$ achieves a maximum.  We may assume that $(x_0, t_0)$ lies in $B_{r_{m+1}} \times (0,T]$, otherwise, by the inductive hypothesis, we are finished.     Suppose first that $f^{m+2} |\nabla^m\Rm|^2 \le  A^{m+2}$ at the point $(x_0,t_0)$.  Then  at any $(x,t) \in \ov{B_{r_{m+2}}} \times [0,T]$ we have 
\begin{equation}
	f^{m+2} |\nabla^m\Rm|^2 
	\myle  A^{m+2} + f^{m+1} B |\nabla^{m-1}\Rm|^2 \Big|_{(x_0,t_0)}, 
\label{maxprin higherorder 1}
\end{equation}
and our claim follows by the inductive hypothesis.  Otherwise we have
\begin{equation}
	f^{m+2} |\nabla^m\Rm|^2 >  A^{m+2} \textrm{ at } (x_0,t_0).
\label{assumed min for grad^m Rm}
\end{equation}

We note that by the inductive hypothesis we always have 
\begin{equation}
	|\nabla^i\Rm| |\nabla^j\Rm| \myle C (A/f)^{\frac{i+j}2+2} \textrm{ when } i,j <m .
\end{equation}
At $(x_0, t_0)$, 
\begin{align}
	0 \le (\partial_t-\Delta) F 
	\myle& C \psi^2 f^{m+1} |\nabla^m\Rm|^2 
		+ |\Delta \psi^2| f^{m+2} | \nabla^m \Rm|^2
		- \psi^2 f^{m+2} |\nabla^{m+1}\Rm|^2   
		\nonumber\\
	& + C \psi^2 f^{m+2} |\Rm| |\nabla^m\Rm|^2 
		+ C \psi^2 f^{m+2} (A/f)^{\frac m 2 + 2} |\nabla^m\Rm|
		\nonumber\\
	& + C f^{m+2} \psi\left|\nabla\psi\right| 
		\left|\nabla^{m+1}\Rm\right| \left|\nabla^m\Rm\right|
		\nonumber\\ 
	& + C B f^{m} (A/f)^{m+1} 
		- B f^{m+1} |\nabla^m \Rm|^2 
		\nonumber\\
	& + C B f^{m+1} |\Rm| |\nabla^{m-1}\Rm|^2
		+ C B f^{m+1} (A/f)^{\frac m 2 + \frac32} |\nabla^{m-1} \Rm|
		\nonumber\\
	\myle& C f^{m+1} A |\nabla^m \Rm|^2 
		+ C f^{\frac{m}{2}} A^{\frac{m}{2}+2} |\nabla^m\Rm| 
		\nonumber\\
	& - C_1 A f^{m+1} |\nabla^m\Rm|^2 
		+ C A^{m+3}f^{-1}
		\nonumber\\
	\myle& - f^{m+1} A |\nabla^m\Rm|^2 + C' A^{m+3}f^{-1}   
\end{align}
where the final inequality follows from (\ref{assumed min for grad^m Rm}) and by taking the uniform constant $C_1$ in $B = C_1 A$ uniformly large enough.  Hence  $f^{m+2} |\nabla^m\Rm|^2 \le C' A^{m+2}$ at $(x_0,t_0)$ and then, arguing in a similar way to (\ref{maxprin higherorder 1}) above, this completes the inductive step.  Thus (\ref{Bound on grad^m Rm}) is established.

\section{Proof of Corollary \ref{corollary}} \label{secproof}

There are various ways to deduce Corollary \ref{corollary} from Theorem \ref{maintheorem}.  
We could directly apply standard local parabolic theory (as discussed in  \cite{Chau, PSSW} for example), or  the method in \cite{CK}.  However, in our setting, we do not even need that $g(t)$ is a solution of a parabolic equation and instead we use an argument similar to one in 
 \cite{SW} which uses only standard linear elliptic theory and some embedding theorems.

Fix a time $t \in (0, T]$.  Regarding $g_{i \bar{j}}$ as a set of $n^2$ functions, we consider the equations
\begin{equation} \label{poisson}
\hat{\Delta} g_{i \bar{j}} = - \sum_k R_{k \bar{k} i \bar{j}} + \sum_{k,p,q} g^{q \bar{p}} \partial_k g_{i \bar{q}} \partial_{\bar{k}} g_{p \bar{j}}=: Q_{i \bar{j}}.
\end{equation}
where $\hat{\Delta} = \sum_k \partial_k \partial_{\bar{k}}$.  For each fixed $i, j$, we can regard (\ref{poisson}) as Poisson's equation $\hat{\Delta} g_{i \bar{j}} = Q_{i \bar{j}}$.

For the purposes of this section we will say that a quantity $Z$ is \emph{uniformly bounded} if there exist constants $C, \alpha, \beta, \gamma$ depending only on $\hat{\omega}$ and $T$ such that 
$Z \le C N^{\alpha} r^{-\beta} t^{-\gamma}$.  In the case when the constants may depend on $\omega_0$, we insist that $\gamma=0$.

Let $r=r_0 > r_1 > \cdots > r/2$ be as above. Fix $p>2n$.  From what we have proved, each  $\| Q_{i \bar{j}} \|_{L^p(B_{r_2})}$ is uniformly bounded.  Applying the standard elliptic estimates for the Poisson equation (see for example Theorem 9.11 of \cite{GT}) to (\ref{poisson}) we see that the Sobolev norm $\| g_{i \bar{j}} \|_{L^p_2(B_{r_3})}$ is uniformly bounded.   Morrey's embedding theorem (Theorem 7.17 of \cite{GT}) gives that  $\| g_{i \bar{j}} \|_{C^{1+ \kappa}(B_{r_4})}$ is uniformly bounded for some $0<\kappa<1$. 

The key observation we now need is that  the $m$th derivative of $Q_{i \bar{j}}$  can be written as a finite sum 
$\sum_{s} A_{s}*B_s$
where each $A_s$ or $B_s$ is either a  covariant derivative of $\textrm{Rm}$  or a quantity involving derivatives of $g$ up to order at most $m+1$.  Hence if $g$ is uniformly bounded in $C^{m+1+ \kappa}$  then each $Q_{i \bar{j}}$ is uniformly bounded in $C^{m+ \kappa}$, after possibly passing to a slightly smaller ball.

Applying this observation with $m=0$ we see that each $\| Q_{i\bar{j}} \|_{C^{\kappa}(B_{r_4})}$ is uniformly bounded. 
The standard Schauder estimates for the Poisson equation (Theorem 4.8 of \cite{GT}) give that $\| g_{i \bar{j}} \|_{C^{2+\kappa}(B_{r_5})}$ is uniformly bounded.  

We can now apply a bootstrapping argument.  Applying the observation with $m=1$ we see that $Q_{i\bar{j}}$  is uniformly bounded in $C^{1+\kappa}$ on a slightly smaller ball and so on. This completes the proof of the corollary.

\bigskip
\noindent
$*$ Department of Mathematics, California Polytechnic State University, San Luis Obispo, CA 93407

\noindent
$\dagger$ Department of Mathematics, University of California San Diego, La Jolla, CA 92093

\end{document}